\documentclass[a4paper,twoside,11pt]{article}
\usepackage{amsfonts}
\usepackage{amssymb}\usepackage{amsmath}
\newtheorem{thm}{Theorem}[section]
 \newtheorem{cor}{Corollary}[section]
 \newtheorem{lem}[thm]{Lemma}

 \numberwithin{equation}{section}
\newcommand{\double}{\baselineskip 1.24 \baselineskip}

\title{A generalized Young inequality and some new results on fractal
space
}

\author{{Guang-Sheng  Chen\thanks{\text{E-mail address}: cgswavelets@126.com(Chen)
}\quad}\\
{\small Department of Computer Engineering, Guangxi Modern
Vocational Technology College,} \\{\small Hechi,Guangxi, 547000,
P.R. China}
}
\begin{document}
\date{}
\maketitle \double

\textbf{Abstract:}\quad Starting with real line number system based on the theory
of the Yang's fractional set, the generalized Young inequality is
established. By using it some results on the generalized inequality in
fractal space are investigated in detail.   \\
\textbf{Keywords:} Dfractal, real line number system, fractional set,
generalized Young inequality\\
\textbf{MSC2010: } 28A80, 26D15

\section{Introduction}
\hskip\parindent
The classical Young inequality [1--4] is not only interesting in itself but
also very useful. The purpose of this work is to establish a generalized
Young inequality on fractional set and other inequality based on it. Start
with, we review the Yang's fractional set and Yang's geometric expressions
for the real line number system.
\subsection{ Theory of the Yang's fractional set}
\hskip\parindent
Recently, the theory of Yang's fractional sets of element sets [6-12] was
introduced as follows:

For $0 < \alpha \le 1$, we have the following $\alpha $-type set of element
sets:

${\rm N}_0^\alpha $ The $\alpha $-type set of the natural numbers are defined as
the set $\{0^\alpha ,1^\alpha ,2^\alpha ,\ldots ,n^\alpha ,\ldots \}$;

${\rm N}_ + ^\alpha $ The $\alpha $-type set of the natural numbers are defined
as the set $\{1^\alpha ,2^\alpha ,\ldots ,n^\alpha ,\ldots \}$;

${\rm Z}_0^\alpha $ The $\alpha $-type set of the integers are defined as the
set $\{0^\alpha ,\pm 1^\alpha ,\pm 2^\alpha ,\ldots ,\pm n^\alpha ,\ldots
\}$;

${\rm Z}^\alpha $ The $\alpha $-type set of the integers are defined as the set
$\{\pm 1^\alpha ,\pm 2^\alpha ,\ldots ,\pm n^\alpha ,\ldots \}$;

${\rm Q}^\alpha $ The $\alpha $-type set of the rational numbers are defined as
the set $\{m^\alpha = (p / q)^\alpha :p,q \in {\rm Z},q \ne 0\}$;

$\Im ^\alpha $ The $\alpha $-type set of the irrational numbers are defined as
the set $\{m^\alpha \ne (p / q)^\alpha :p,q \in {\rm Z},q \ne 0\}$;

${\rm R}^\alpha $ The $\alpha $-type set of the real line numbers are defined as
the set ${\rm R}^\alpha = {\rm Q}^\alpha \cup \Im ^\alpha $.

\subsection{Yang's geometric expressions for the real line number system}
\hskip\parindent
Geometric representation of real line numbers on a fractional set as points
on a real line called the real line axis. For each real line number there
correspond one and only one point on the real line[6].

For example, $1^\alpha + 2^\alpha = 3^\alpha $. That is, the geometric
representation is that cantor set $[0,3]$ is equivalent to the sum of cantor
set $[0,1]$ and cantor set $[1,3]$. The dimension of cantor set is $\alpha $
, for $0 < \alpha \le 1$ and, $1^\alpha $ , $2^\alpha $ and $3^\alpha $ are
real line numbers on a fractional set. If $a^\alpha ,b^\alpha ,c^\alpha $
belong to the set ${\rm R}^\alpha $ of real line numbers, then we have the
following operation:

(1) $a^\alpha + b^\alpha $ and $a^\alpha b^\alpha $ belong to the set ${\rm
R}^\alpha $

(2) $a^\alpha + b^\alpha = b^\alpha + a^\alpha = (a + b)^\alpha = (b +
a)^\alpha $;

(3) $a^\alpha + (b^\alpha + c^\alpha ) = (a^\alpha + b^\alpha ) + c^\alpha $;

(4) $a^\alpha b^\alpha = b^\alpha a^\alpha = (ab)^\alpha = (ba)^\alpha $;

(5) $a^\alpha (b^\alpha c^\alpha ) = (a^\alpha b^\alpha )c^\alpha $;

(6) $a^\alpha (b^\alpha + c^\alpha ) = a^\alpha b^\alpha + a^\alpha c^\alpha
$;

(7) $a^\alpha + 0^\alpha = 0^\alpha + a^\alpha = a^\alpha $ and $a^\alpha
\cdot 1^\alpha = 1^\alpha \cdot a^\alpha = a^\alpha $ .

If $a^\alpha - b^\alpha $ is a nonnegative number, we say that $a^\alpha $
is greater than or equal to $b^\alpha $ or $b^\alpha $ is less than or equal
to $a^\alpha $ , and write, respectively, $a^\alpha \ge b^\alpha $ or
$b^\alpha \le a^\alpha $ . If there is no possibility that $a^\alpha =
b^\alpha $ , we write

$a^\alpha > b^\alpha $ or $b^\alpha < a^\alpha $ .

Geometrically, $a^\alpha > b^\alpha $ if the point on the real line axis
corresponding to $a^\alpha $ lies to the left of the point corresponding to
$b^\alpha $ .

Suppose that $a^\alpha $ , $b^\alpha $ and $c^\alpha $ are any given real
line numbers, then we have the following relations:

(1) Either $a^\alpha > b^\alpha $ , $a^\alpha = b^\alpha $ or $a^\alpha <
b^\alpha $ ( Law of trichotomy );

(2) If $a^\alpha > b^\alpha $ and $b^\alpha > c^\alpha $ , then $a^\alpha >
c^\alpha $ (Law of transitivity);

(3) If $a^\alpha > b^\alpha $ , then $a^\alpha + c^\alpha > b^\alpha +
c^\alpha $ ;

(4) If $a^\alpha > b^\alpha $ and $c^\alpha > 0^\alpha $ , then $a^\alpha
c^\alpha > b^\alpha c^\alpha $ ;

(5) If $a^\alpha > b^\alpha $ and $c^\alpha < 0^\alpha $ , then $a^\alpha
c^\alpha < b^\alpha c^\alpha $ .

The formula is similar to classical one in case of $\alpha = 1$. As direct
results, the following inequalities are valid:

If $a^\alpha > b^\alpha $ , then $a > b$;

If $a^\alpha = b^\alpha $ , then $a = b$;

If $a^\alpha < b^\alpha $ , then $a < b$.

\section{The generalized Young inequality }
\hskip\parindent
In the section, we give the proof of the generalized Young inequality. Here
we first start with the generalized Bernoulli's inequality.
\begin{lem}\label{Lemma 2.1} (Generalized Bernoulli's inequality) Let $y > 0$，

(1) for $0 < m < 1$, then
\begin{equation*}
y^{\alpha m} - 1^\alpha \le m(y - 1)^\alpha .
\tag{2.1}
\label{2.1}
\end{equation*}

(2) for $m > 1$, then
\begin{equation*}
y^{\alpha m} - 1^\alpha \ge m(y - 1)^\alpha  .
\tag{2.2}
\label{2.2}
\end{equation*}
\end{lem}
\textbf{ Remark 1.} This is classical Bernoulli's inequality in case of fractal
dimension $\alpha = 1$[5].
\begin{thm}\label{Theorem 2.2}. Let $a$, $b \ge 0$, $p > 1$, $\frac{1}{p} + \frac{1}{q}
= 1$,then
\begin{equation*}
a^\alpha b^\alpha \le \frac{a^{p\alpha }}{p} + \frac{b^{q\alpha }}{q},
\tag{2.3}
\label{2.3}
\end{equation*}
which is equality holding if and only if $a^p = b^q$.
\end{thm}
\textbf{Proof}. Setting $y = a / b$, by (2.1) we have
\begin{equation*}
(a / b)^{\alpha m} - 1^\alpha \le m(a / b - 1)^\alpha .
\tag{2.4}
\label{2.4}
\end{equation*}
Multiplying both sides by $b^\alpha $ in (2.4) gives
\begin{equation*}
a^{\alpha m}b^{\alpha (1 - m)} - b^\alpha \le m(a - b)^\alpha .
\tag{2.5}
\label{2.5}
\end{equation*}
Then, we obtain that
\begin{equation*}
a^{\alpha m}b^{\alpha (1 - m)} \le ma^\alpha + (1 - m)b^\alpha .
\tag{2.6}
\label{2.6}
\end{equation*}
Let $m = \frac{1}{p} < 1$, then we directly deduce (2.3).
\\ \textbf{ Remark 2.} (2.3) was proposed[6], here we give its proof.it is classical Young inequality in case of fractal dimension $\alpha=1$.

\begin{thm}\label{Theorem 2.3.} Let $a$,$b \ge 0$, $0 < p < 1$ and $\frac{1}{p} +
\frac{1}{q} = 1$, then
\begin{equation*}
a^\alpha b^\alpha \ge \frac{a^{p\alpha }}{p} + \frac{b^{q\alpha }}{q},
\tag{2.7}
\label{2.7}
\end{equation*}
where equality holds if and only if $a^p = b^q$.
\end{thm}
\textbf{Proof}. For $x$, $y \ge 0$, $0 < p < 1$, by (2.3), we have
\begin{equation*}
x^\alpha y^\alpha \le px^{\textstyle{\alpha \over p}} + (1 -
p)y^{\textstyle{\alpha \over q}}.
\tag{2.8}
\label{2.8}
\end{equation*}
Set $a = \frac{1}{p^p}x^py^p$, $b = \frac{1}{p^{1 - p}}y^{ - p}$. By (2.8) we
obtain
\[
a^\alpha b^\alpha \ge \frac{a^{p\alpha }}{p} + \frac{b^{q\alpha }}{q}.
\]
Hence we complete the proof of Theorem 2.2.

As a direct result, we have the following:
\begin{cor}\label{corollary  2.1} Let $a_i \ge 0$, $p_i \in {\rm R}$, $i = 1,2,\ldots
n$, $\sum\limits_{i = 1}^n {1 / p_i } = 1$,

(1) for $p_i > 1$, we have
\begin{equation*}
\prod\limits_{i = 1}^n {a_i^\alpha } \le \sum\limits_{i = 1}^n {a_i^{p_i
\alpha } / p_i },
\tag{2.9}
\label{2.9}
\end{equation*}
where the equality holds if $a_j^{p_j} = a_k^{p_k}$ for $\forall j,k$ .

(2) for $0 < p_1 < 1$, $p_i < 0$, $i = 2,\ldots n$,we have
\begin{equation*}
\prod\limits_{i = 1}^n {a_i^\alpha } \ge \sum\limits_{i = 1}^n {a_i^{p_i
\alpha } / p_i },
\tag{2.10}
\label{2.10}
\end{equation*}
where the equality holds if $a_j^{p_j} = a_k^{p_k}$ for $\forall j,k$ .
\end{cor}
\section{Useful results based on generalized Young inequality}
\hskip\parindent
In the section we discuss some generalizations of holder inequality and
Minkowski inequality. In order to proof our results, we first review the
h\"{o}lder inequality and Minkowski inequality [6]:
\begin{thm}\label{Theorem 3.1}.(see [6]) Let $\left| {x_i } \right|$, $\left| {y_i } \right|
\ge 0$, $p > 1$, $\frac{1}{p} + \frac{1}{q} = 1$, $i = 1,2,\ldots n$, then
\begin{equation*}
\sum\limits_{i = 1}^n {\left| {x_i^\alpha } \right|\left| {y_i^\alpha }
\right|} \le \left( {\sum\limits_{i = 1}^n {\left| {x_i^\alpha } \right|^p}
} \right)^{1 / p}\left( {\sum\limits_{i = 1}^n {\left| {y_i^\alpha }
\right|^q} } \right)^{1 / q}.
\tag{3.1}
\label{3.1}
\end{equation*}
Equalities holding if and only if $\lambda_1\left| {x_i } \right|=\lambda_2\left| {y_i } \right|$, where $\lambda_1$ and $\lambda_2$ are constants.
\end{thm}
\begin{thm}\label{Theorem 3.2}.(see[6]) Let $\left| {x_i } \right|$, $\left| {y_i } \right|
\ge 0$, $p > 1$, $i = 1,2,\ldots n$, then
\begin{equation*}
\left( {\sum\limits_{i = 1}^n {\left| {x_i^\alpha + y_i^\alpha } \right|^p}
} \right)^{1 / p} \le \left( {\sum\limits_{i = 1}^n {\left| {x_i^\alpha }
\right|^p} } \right)^{1 / p} + \left( {\sum\limits_{i = 1}^n {\left|
{y_i^\alpha } \right|^q} } \right)^{1 / q}.
\tag{3.2}
\label{3.2}
\end{equation*}
Equalities holding if and only if $\lambda_1\left| {x_i } \right|=\lambda_2\left| {y_i } \right|$, where $\lambda_1$ and $\lambda_2$ are constants.
\end{thm}
\begin{thm}\label{Theorem 3.3}. Let $\left| {x_i } \right|$, $\left| {y_i } \right| \ge
0$, $0 < p < 1$, $\frac{1}{p} + \frac{1}{q} = 1$, $i = 1,2,\ldots n$, then
\begin{equation*}
\sum\limits_{i = 1}^n {\left| {x_i^\alpha } \right|\left| {y_i^\alpha }
\right|} \ge \left( {\sum\limits_{i = 1}^n {\left| {x_i^\alpha } \right|^p}
} \right)^{1 / p}\left( {\sum\limits_{i = 1}^n {\left| {y_i^\alpha }
\right|^q} } \right)^{1 / q}.
\tag{3.3}
\label{3.3}
\end{equation*}
Equalities holding if and only if $\lambda_1\left| {x_i } \right|=\lambda_2\left| {y_i } \right|$, where $\lambda_1$ and $\lambda_2$ are constants.
\end{thm}
\textbf{Proof. }Set $c = 1 / p$, then we have $q = - pd$, $d = c / (c - 1)$.
By (3.1), we obtain
\begin{equation*}
\begin{split}
 &\sum\limits_{i = 1}^n {\left| {x_i^\alpha } \right|^p} = \sum\limits_{i =
1}^n {\left| {x_i^\alpha y_i^\alpha } \right|^p\left| {y_i^\alpha }
\right|^{ - p}} \\
& \le \left( {\sum\limits_{i = 1}^n {\left| {x_i^\alpha y_i^\alpha }
\right|^{pc}} } \right)^{1 / c}\left( {\sum\limits_{i = 1}^n {\left|
{y_i^\alpha } \right|^{ - pd}} } \right)^{1 / d} \\
 &= \left( {\sum\limits_{i = 1}^n {\left| {x_i^\alpha y_i^\alpha } \right|} }
\right)^{1 / p}\left( {\sum\limits_{i = 1}^n {\left| {y_i^\alpha }
\right|^q} } \right)^{1 - p}. \\
\end{split}
\tag{3.4}
\label{3.4}
\end{equation*}
In (3.4), multiplying both sides by $\left( {\sum\limits_{i = 1}^n {\left|
{y_i^\alpha } \right|^q} } \right)^{p - 1}$ yields
\begin{equation*}
\sum\limits_{i = 1}^n {\left| {x_i^\alpha } \right|^p} \left(
{\sum\limits_{i = 1}^n {\left| {y_i^\alpha } \right|^q} } \right)^{p - 1}
\le \left( {\sum\limits_{i = 1}^n {\left| {x_i^\alpha y_i^\alpha } \right|}
} \right)^{1 / p}.
\tag{3.5}
\label{3.5}
\end{equation*}
Using (3.5) implies that
\[
\sum\limits_{i = 1}^n {\left| {x_i^\alpha } \right|\left| {y_i^\alpha }
\right|} \ge \left( {\sum\limits_{i = 1}^n {\left| {x_i^\alpha } \right|^p}
} \right)^{1 / p}\left( {\sum\limits_{i = 1}^n {\left| {y_i^\alpha }
\right|^q} } \right)^{1 / q}.
\]
\begin{thm}\label{Theorem 3.4}. Let $\left| {x_i } \right|$,$\left| {y_i } \right| \ge
0$, $0 < p < 1$,$i = 1,2,\ldots n$, then
\begin{equation*}
\left( {\sum\limits_{i = 1}^n {\left| {x_i^\alpha + y_i^\alpha } \right|^p}
} \right)^{1 / p} \ge \left( {\sum\limits_{i = 1}^n {\left| {x_i^\alpha }
\right|^p} } \right)^{1 / p} + \left( {\sum\limits_{i = 1}^n {\left|
{y_i^\alpha } \right|^p} } \right)^{1 / p}.
\tag{3.6}
\label{3.6}
\end{equation*}
Equalities holding if and only if $\lambda_1\left| {x_i } \right|=\lambda_2\left| {y_i } \right|$, where $\lambda_1$ and $\lambda_2$ are constants.
\end{thm}
\textbf{Proof.} $A_n = \sum\limits_{i = 1}^n {\left| {x_i^\alpha }
\right|^p} $，$B_n = \sum\limits_{i = 1}^n {\left| {y_i^\alpha } \right|^p}
$,
$
C_n = \left( {\sum\limits_{i = 1}^n {\left| {x_i^\alpha } \right|^p} }
\right)^{1 / p} + \left( {\sum\limits_{i = 1}^n {\left| {y_i^\alpha }
\right|^q} } \right)^{1 / q}
$, by H\"{o}lder inequality, in view of $0 < p < 1$, we have
\begin{equation*}
\begin{split}
 &C_n = \sum\limits_{i = 1}^n {(\left| {x_i^\alpha } \right|^pA_n^{1 / p - 1}
+ \left| {y_i^\alpha } \right|^pB_n^{1 / p - 1} )} \\
 &\le \sum\limits_{i = 1}^n {\left| {x_i^\alpha + y_i^\alpha }
\right|^p(A_n^{1 / p} + B_n^{1 / p} )^{1 - p}} = C_n^{1 - p} \sum\limits_{i
= 1}^n {\left| {x_i^\alpha + y_i^\alpha } \right|^p} .\\
 \end{split}
\tag{3.7}
\label{3.7}
\end{equation*}
By inequality (3.7), we arrive at reverse Minkowski's inequality and the
theorem is completely proved.
\begin{cor}\label{Collory 3.5} Let $\left| {x_{ij} } \right| \ge 0$, $p_j \in {\rm
R}$, $i = 1,2,\ldots n$, $j = 1,2,\ldots m\sum\limits_{j = 1}^m {1 / p_j } =
1$.

(1) for $p_j > 1$, we have
\begin{equation*}
\sum\limits_{i = 1}^n {\prod\limits_{j = 1}^m {\left| {x_{ij}^\alpha }
\right|} } \le \prod\limits_{j = 1}^m {\left( {\sum\limits_{i = 1}^n {\left|
{x_{ij}^\alpha } \right|^{p_j }} } \right)^{1 / p_j }}.
\tag{3.8}
\label{3.8}
\end{equation*}
Equalities holding if and only if $\lambda_j\left| {x_{ij }} \right|=\lambda_k\left| {x_{ik} } \right|$ for $\forall j,k$, where $\lambda_j$ and $\lambda_k$ are constants.

(2) for $0 < p_1 < 1$, $p_j < 0$, $j = 2,\ldots m$,we have
\begin{equation*}
\sum\limits_{i = 1}^n {\prod\limits_{j = 1}^m {\left| {x_{ij}^\alpha }
\right|} } \ge \prod\limits_{j = 1}^m {\left( {\sum\limits_{i = 1}^n {\left|
{x_{ij}^\alpha } \right|^{p_j }} } \right)^{1 / p_j }}.
\tag{3.9}
\label{3.9}
\end{equation*}
\end{cor}
Equalities holding if and only if $\lambda_j\left| {x_{ij} } \right|=\lambda_k\left| {x_{ik} } \right|$ for $\forall j,k$, where $\lambda_j$ and $\lambda_k$ are constants.

\begin{cor}\label{Collory 3.6} Let $\left| {x_{ij} } \right| \ge 0$, $i = 1,2,\ldots
n$,$j = 1,2,\ldots m$，,

(1)for $p > 1$, we have

\begin{equation*}
\left( {\sum\limits_{i = 1}^n {\left| {\sum\limits_{j = 1}^m {x_{ij}^\alpha
} } \right|^p} } \right)^{1 / p} \le \sum\limits_{i = 1}^n {\left(
{\sum\limits_{j = 1}^m {\left| {x_{ij}^\alpha } \right|^p} } \right)^{1 /
p}}.
\tag{3.10}
\label{3.10}
\end{equation*}
Equalities holding if and only if $\lambda_j\left| {x_{ij} } \right|=\lambda_k\left| {x_{ik} } \right|$ for $\forall j,k$, where $\lambda_j$ and $\lambda_k$ are constants.

(2)for $0 < p < 1$ ,we have
\begin{equation*}
\left( {\sum\limits_{i = 1}^n {\left| {\sum\limits_{j = 1}^m {x_{ij}^\alpha
} } \right|^p} } \right)^{1 / p} \ge \sum\limits_{i = 1}^n {\left(
{\sum\limits_{j = 1}^m {\left| {x_{ij}^\alpha } \right|^p} } \right)^{1 /
p}}.
\tag{3.11}
\label{3.11}
\end{equation*}
Equalities holding if and only if $\lambda_j\left| {x_{ij} } \right|=\lambda_k\left| {x_{ik} } \right|$ for $\forall j,k$, where $\lambda_j$ and $\lambda_k$ are constants.
\end{cor}
\begin{thm}\label{Theorem 3.7} Let $\left| {x_i } \right|$, $\left| {y_i } \right| \ge
0$ , $0 < r < 1 < p$, $i = 1,2,\ldots n$, then
\begin{equation*}
\left( {\frac{\sum\limits_{i = 1}^n {\left| {x_i^\alpha + y_i^\alpha }
\right|^p} }{\sum\limits_{i = 1}^n {\left| {x_i^\alpha + y_i^\alpha }
\right|^r} }} \right)^{1 / (p - r)} \le \left( {\frac{\sum\limits_{i = 1}^n
{\left| {x_i^\alpha } \right|^p} }{\sum\limits_{i = 1}^n {\left| {x_i^\alpha
} \right|^r} }} \right)^{1 / (p - r)} + \left( {\frac{\sum\limits_{i = 1}^n
{\left| {y_i^\alpha } \right|^p} }{\sum\limits_{i = 1}^n {\left| {y_i^\alpha
} \right|^r} }} \right)^{1 / (p - r)}.
\tag{3.12}
\label{3.12}
\end{equation*}
Equalities holding if and only if $\lambda_1\left| {x_i} \right|=\lambda_2\left|{y_i}\right|$,where $\lambda_1$ and $\lambda_2$ are constants.
\end{thm}
\textbf{Proof. }By Theorem 3.1 and Theorem 3.2， We have
\begin{equation*}
\begin{split}
& \left( {\sum\limits_{i = 1}^n {\left| {x_i^\alpha + y_i^\alpha } \right|^p}
} \right)^{1 / (p - r)} \le \left( {\left( {\sum\limits_{i = 1}^n {\left|
{x_i^\alpha } \right|^p} } \right)^{1 / p} + \left( {\sum\limits_{i = 1}^n
{\left| {y_i^\alpha } \right|^p} } \right)^{1 / p}} \right)^{p / (p - r)} \\
 &= \left( {\left( {\frac{\sum\limits_{i = 1}^n {\left| {x_i^\alpha }
\right|^p} }{\sum\limits_{i = 1}^n {\left| {x_i^\alpha } \right|^r} }}
\right)^{1 / p}\left( {\sum\limits_{i = 1}^n {\left| {x_i^\alpha }
\right|^r} } \right)^{1 / p} + \left( {\frac{\sum\limits_{i = 1}^n {\left|
{y_i^\alpha } \right|^p} }{\sum\limits_{i = 1}^n {\left| {y_i^\alpha }
\right|^r} }} \right)^{1 / p}\left( {\sum\limits_{i = 1}^n {\left|
{y_i^\alpha } \right|^r} } \right)^{1 / p}} \right)^{p / (p - r)} \\
 &\le \left( {\left( {\frac{\sum\limits_{i = 1}^n {\left| {x_i^\alpha }
\right|^p} }{\sum\limits_{i = 1}^n {\left| {x_i^\alpha } \right|^r} }}
\right)^{1 / (p - r)} + \left( {\frac{\sum\limits_{i = 1}^n {\left|
{y_i^\alpha } \right|^p} }{\sum\limits_{i = 1}^n {\left| {y_i^\alpha }
\right|^r} }} \right)^{1 / (p - r)}} \right)\left( {\left( {\sum\limits_{i =
1}^n {\left| {x_i^\alpha } \right|^r} } \right)^{1 / r} + \left(
{\sum\limits_{i = 1}^n {\left| {y_i^\alpha } \right|^r} } \right)^{1 / r}}
\right)^{r / (p - r)}. \\
 \end{split}
\tag{3.13}
\label{3.13}
\end{equation*}
Using reverse Minkowski inequality implies that
\begin{equation*}
\left( {\left( {\sum\limits_{i = 1}^n {\left| {x_i^\alpha } \right|^r} }
\right)^{1 / r} + \left( {\sum\limits_{i = 1}^n {\left| {y_i^\alpha }
\right|^r} } \right)^{1 / r}} \right)^r \le \sum\limits_{i = 1}^n {\left|
{x_i^\alpha + y_i^\alpha } \right|^r}.
\tag{3.14}
\label{3.14}
\end{equation*}
By (3.13) and (3.14), we get (3.12). Hence, the theorem is completely proved.

\begin{cor}\label{Collory 3.8} Let $\left| {x_{ij} } \right| \ge 0$, $0 < r < 1 <
p$, $i = 1,2,\ldots n$,$j = 1,2,\ldots m$ then
\begin{equation*}
\left( {\frac{\sum\limits_{i = 1}^n {\left| {\sum\limits_{j = 1}^m
{x_{ij}^\alpha } } \right|^p} }{\sum\limits_{i = 1}^n {\left|
{\sum\limits_{j = 1}^m {x_{ij}^\alpha } } \right|^r} }} \right)^{1 / (p -
r)} \le \sum\limits_{i = 1}^n {\left( {\frac{\sum\limits_{j = 1}^m {\left|
{x_{ij}^\alpha } \right|^p} }{\sum\limits_{j = 1}^m {\left| {x_{ij}^\alpha }
\right|^r} }} \right)^{1 / (p - r)}}.
\tag{3.15}
\label{3.15}
\end{equation*}
Equalities holding if and only if $\lambda_j\left| {x_{ij} } \right|=\lambda_k\left| {x_{ik} } \right|$ for $\forall j,k$, where $\lambda_j$ and $\lambda_k$ are constants.
\end{cor}


\begin{thebibliography}{99}\small

\bibitem{1} O. H\"{o}lder, Uber einen Mittelwerthssatz, Nachr. Ges. Wiss. Gottingen (1889) 38-47.
\bibitem{2} D.S. Mitrinovic, Analytic Inequalities, Springer-Verlag, New York, Heidelberg, Berlin, 1994.
\bibitem{3}T. Takahashi, Remarks on some inequalities, T\^{o}hoku Math. J. 36 (1932) 99-106.
\bibitem{4}W.H. Young, On class of summable functions and there Fourier series, Proc. Roy. Soc. London A 87 (1912) 225-229.
\bibitem{5} J. Kuang. Applied Inequalities, Shandong Science Press, Jinan,
2003.

\bibitem{6}X. Yang, Local Fractional Integral Transforms, Progress in Nonlinear Science, 4(2011): 1-225

\bibitem{7} X.Yang, Local fractional Laplace's transform based on the local fractional calculus. In: Proc. of The 2011 International Conference on Computer Science and Information Engineering (CSIE2011), Springer, 2011.

\bibitem{8}X.Yang, Z.Kang, C.Liu. Local fractional Fourier's transform based on the local fractional calculus, In: Proc. of The International Conference on Electrical and Control Engineering (ICECE 2010), 2010,1242-1245.

\bibitem{9} F. Gao, X.Yang, Z. Kang. Local fractional Newton's method derived from modified local fractional calculus. In: Proc. of The Second Scientific and Engineering Computing Symposium on
Computational Sciences and Optimization (CSO 2009), 2009, 228-232.

\bibitem{10}X.Yang, F. Gao. The Fundamentals of local fractional derivative of the  one-variable non-differentiable functions, World Sci-Tech R\&D, 31(5), 2009, 920-921.

\bibitem{11} X.Yang, L.Li, R.Yang. Problems of local fractional definite integral of the one-variable
non-differentiable function, World Sci-Tech R\&D, 31(4), 2009, 722-724.

\bibitem{12} X.Yang, F.Gao. Fundamentals of local fractional iteration of the continuously nondifferentiable functions derived from local fractional calculus. In: Proc. of The 2011 International Conference on Computer Science and Information Engineering (CSIE2011), Springer, 2011.
\end{thebibliography}
 \end{document}